\documentclass[a4paper,10pt]{article}

\usepackage{amsmath,amssymb,amsthm}
\usepackage{graphicx}
\usepackage{algorithm,algorithmic}
\usepackage{tikz,pgfpages}

\usepackage[top=1in, bottom=1in, left=1.25in, right=1.25in]{geometry}

\newcommand{\normtwo}[1]{\lVert #1 \rVert_2}
\newcommand{\normB}[1]{\lVert #1 \rVert_B}
\newcommand{\define}{\stackrel{\text{def}}{=}}
\newcommand{\prior}{\Gamma_\text{sys}}
\newcommand{\noise}{\Gamma_\text{noise}}
\newcommand{\bigO}{{\cal{O}}}

\newcommand{\nlogn}{N_s\log^\beta N_s}
\newcommand{\bx}{\textbf{x}}
\newcommand{\at}[2]{\left.#1\right|_{#2}}
\newcommand{\normal}{ {\cal{N}}}

\newcommand{\logdet}{\log\text{det}  }
\newcommand{\trace}{\text{Trace}}

\newcommand{\Span}[1]{\text{span}\left\{#1\right\}}

\newcommand{\skk}{\hat{s}_{k|k}}
\newcommand{\skkp}{\hat{s}_{k+1|k}}
\newcommand{\skkpp}{\hat{s}_{k+1|k+1}}

\newcommand{\Skk}{\Sigma_{k|k}}
\newcommand{\Skkp}{\Sigma_{k+1|k}}
\newcommand{\Skkpp}{\Sigma_{k+1|k+1}}

\newcommand*{\xMin}{0}%
\newcommand*{\xMax}{10}%
\newcommand*{\yMin}{0}%
\newcommand*{\yMax}{10}%

\title{Fast Kalman Filter  using Hierarchical-matrices and low-rank perturbative approach}
\author{Arvind K. Saibaba \thanks{  Department of Electrical and Computer Engineering, Tufts University } \and Eric L. Miller \footnotemark[1] \and  Peter K. Kitanidis \thanks{Institute for Computational and Mathematical Engineering and Department of Civil and Environmental Engineering, Stanford University} }

\begin{document}
\maketitle

\begin{abstract}
We develop a fast algorithm for Kalman Filter applied to the random walk forecast model. The key idea is an efficient representation of the estimate covariance matrix at each time-step as a weighted sum of two contributions - the process noise covariance matrix and a low rank term computed from a generalized eigenvalue problem, which combines information from the noise covariance matrix and the data. We describe an efficient algorithm to update the weights of the above terms and the computation of eigenmodes of the generalized eigenvalue problem (GEP). The resulting algorithm for the Kalman filter with a random walk forecast model scales as $\bigO(N)$ in memory and $\bigO(N \log N)$ in computational cost, where $N$ is the number of grid points. We show how to efficiently compute measures of uncertainty and conditional realizations from the state distribution at each time step. An extension to the case with nonlinear measurement operators is also discussed. Numerical experiments demonstrate the performance of our algorithms, which are applied to a synthetic example from monitoring CO$_2$ in the subsurface using travel time tomography.

\end{abstract}

\section{Introduction}\label{sec:intro}

Kalman filtering is a fundamental tool in statistical time series analysis used to estimate the states of large-scale dynamical systems for which noisy observations are available. In several geophysical and biomedical applications we wish to estimate high-dimensional system parameters. Standard implementations of the Kalman filter are prohibitive because they require $\bigO(N^2)$ in memory and $\bigO(N^3)$ in computational cost, where $N$ is the dimension of the state variable. For large problem sizes, this cost is prohibitively expensive. The main bottleneck for scalable implementation of the Kalman filter is the computation and representation of the state covariance matrix. In this work, we focus our attention on the random walk forecast model. This model is useful in practical applications in which data is acquired at a rapid rate and when changes in states between times when data is obtained are small that they can approximated by a  random walk process. Previous work has considered this model for filtering in the context of time-lapse electrical impedance tomography~\cite{vauhkonen1998kalman,kim2001image,soleimani2007dynamic}, electrical resistivity tomography~\cite{nenna2011application} and CO$_2$ monitoring using seismic travel time tomography~\cite{li2014kalman}.

Several alternatives have been proposed to deal with the computational cost associated with the Kalman filter. One such approach is to represent the state covariance matrices $\Skk$ by a sparse matrix, banded or a strongly tapered structure and update it using sparse matrix techniques. While this approach is effective in certain cases~\cite{furrer2007estimation}, there is no reason to assume  \emph{a priori} that the matrices $\Skk$ will continue to inherit a sparsity structure. A second popular approach is to construct a low-rank approximation to the state covariance matrices using ensemble averages composed of several realizations. This Monte Carlo based approach is known as the Ensemble Kalman filter and is widely used in several areas, in particular, numerical weather prediction~\cite{evensen1994sequential,burgers1998analysis}. Although the ensemble methods provide a cheaper alternative, to get accurate results a large sample size is required which greatly increases the computational costs~\cite{li2014kalman}.

\textbf{Contributions}: In our approach, we consider an efficient representation of the posterior covariance matrix as a low-rank perturbation of the system noise covariance matrix and can be written as $\Skk = \alpha_k\prior - W_kD_kW_k^T$, with $W_k$ chosen such that $W_k^T\prior^{-1} W_k = I$ and $D_k$ is a diagonal matrix. The system noise covariance matrix $\prior$ arising out of Mat\'{e}rn covariance kernels can be efficiently represented using the Hierarchical matrix approach~\cite{saibaba2012efficient,saibaba2012application,ambikasaran2012large} and updates to $\alpha_k$, $D_k$ and $W_k$ are calculated efficiently by solving a generalized eigenvalue problem and repeated application of the Sherman-Morrison-Woodbury update. From this representation, we will show how to compute several measures of uncertainty based on the state covariance matrix. The resulting algorithm for the Kalman filter with Random walk forecast model will be shown to scale as $\bigO(N)$ in memory and $\bigO(N \log N)$ in computational cost. For a small number of measurements, this procedure can be made exact. However, as the number of measurements increase, for several choices of measurement operators and noise covariance matrices, the spectrum of the generalized Hermitian eigenvalue problem decays rapidly and we are justified in only retaining the dominant eigenmodes. We discuss tradeoffs between accuracy and computational cost. The resulting algorithms are applied to a synthetic application to  continuously track CO$_2$ plume in the subsurface using seismic travel times measurements. 

An important contribution of this paper is the development of computationally efficient methods for quantification of uncertainty using optimality measures that have been previously been developed in the context of experimental design, in order to determine the best possible measurements or measurement types, numbers, locations and experimental conditions. Optimal experimental design involves computing measures of conditional uncertainty and minimizes these measures for an optimal design~\cite{alexanderian2013optimal}. Following the approach developed in~\cite{saibaba2013uncertainty} for inverse problems, we extend the computational techniques to computing measures of uncertainty for the data assimilation problem. A related approach has been used to quantify uncertainty in the context of the $4$D-Var data assimilation problem~\cite{singh2013practical}. Conditional realizations from the posterior distribution are also frequently used to understand the uncertainty associated with the state estimation and their computation can be achieved at a cost similar to that of the ``best estimate''. We will demonstrate that conditional realizations computed from the system noise covariance matrix (the most expensive part, which can be computed offline) can be propagated through the time history with only a little additional cost. 

For these reasons, the scalability of our algorithm opens up the possibility of real-time adaptive experimental design and optimal control in systems of much larger dimension than was previously feasible.

\textbf{Related work}: In particular,~\cite{li2014kalman} also consider efficient updates to the Kalman filter using Hierarchical matrices. In order to reduce the memory and computational costs involved with the posterior covariance matrices, they cleverly show how to rewrite the updates of the Kalman filter by only storing and updating the cross-covariance matrices. Further, they show how to compute uncertainty measures such as the variance by computing the diagonals of the posterior covariance matrices. However, their approach cannot be easily extended to compute other measures of uncertainty and non-constant measurement operators. 

Our approach is most similar to the work described in~\cite{paninski2010fast,pnevmatikakis2013fast}, in which the authors represent the posterior covariance matrix as a low-rank perturbation of an equilibrium covariance matrix which is obtained by the solution of a discrete Lyapunov equation. However, a major limitation of their work is that an explicit solution to the equilibrium covariance matrix is only available under special circumstances. In fact, they make the assumption that the system dynamic matrix $F$ is normal and commutes with the system noise covariance matrix $\prior$ and the long term dynamics of covariance matrix (in the absence of observations) converges to an equilibrium covariance matrix. These requirements taken together are restrictive because practical applications of interest do not fully satisfy these criteria. It should be noted that the random walk assumption satisfies the first two requirements (normal and commutes with system noise covariance matrix) but does not have an equilibrium point\footnote{An extension of our algorithms to the case that $F$ is normal and commutes with $\prior$ is straightforward and will not be described.}. As a result we cannot directly apply follow their approach.

The article is organized as follows. In section~\ref{sec:fastkalman} we review the standard computational implementation of the Kalman filter and discuss the computational tools that we will need to derive a fast Kalman filter. Then we derive the fast Kalman filter treating the case of time-invariant measurement operator and then extending our approach to the time-varying measurement operator. In section~\ref{sec:uncert}, we show how to compute several measures of uncertainty based on the {\it a posteriori} state covariance matrix. Furthermore, we show how to compute conditional realizations that can be propagated efficiently in time. Finally, in section~\ref{sec:numerical} we describe the synthetic travel time tomography data assimilation problem and show the performance of our algorithm in comparison to the standard Kalman filter and the Ensemble Kalman Filter. We will demonstrate that the performance of our algorithm has a comparable accuracy to the standard Kalman filter but has a significantly superior computational performance.  

\section{Fast Kalman filter}\label{sec:fastkalman}

\subsection{Problem statement}
Let us begin by reviewing the standard implementation of the  Kalman filter~\cite{kalman1960new}. We denote by $s_k$ and $y_k$ the state variable and observations at step $k$. We assume that $s_k$ and $y_k$ satisfy the following dynamical system 
\begin{align*}
s_{k+1} = & \quad F_ks_{k} + w_k  &\qquad w_k \sim {\cal N} (0,\prior)  \\ \nonumber
y_{k+1} = & \quad H_k s_{k+1} + v_k  &\qquad v_k \sim {\cal N} (0,\noise)  \\ \nonumber
\end{align*}

The system noise $w_k$ is modeled as a Gaussian process with zero mean and covariance $\prior$. $H_k$ is the observation matrix, also called the {\it measurement operator}. The measurements $y_k$ are assumed to be corrupted by noise, which we model as a Gaussian with zero mean and covariance $\noise$.

The Kalman filter~\cite{kalman1960new} is often written out in two stages, 1) prediction stage, in which the state estimate at the previous time step is used to produce an estimate of the state at the current time step, and 2) update stage, in which the prediction is combined with the observation to refine the state estimate. Let $\hat{s}_{k_2|k_1}$ and $\Sigma_{k_2|k_1}$ denote the estimate and covariance (respectively) at step $k_2$ given information till step $k_1$. The equations for the Kalman filter can be summarized as, 

\begin{algorithm}
\begin{algorithmic}
\STATE\textbf{Predict}
\begin{align}
\hat{s}_{k+1|k} \quad = & \quad F_k\hat{s}_{k|k} & - \\  \label{eqn:prediction}
\Skkp \quad = & \quad F_k\Skk F_k^T  + \prior  & \bigO(n_s^3) 
\end{align}

\STATE \textbf{Update}
\begin{align}
S_k \quad = &  \quad H_k\Skkp H_k^T + \noise & \bigO(n_mn_s^2)  \label{eqn:residualcovariance}\\ 
K_k \quad = & \quad \Skkp H^TS_k^{-1} &   \bigO(n_mn_s^2) \label{eqn:kalmangain}\\  
\skkpp \quad = & \quad \skkp+ K_k( y_k - H_k\skkp ) &  \bigO(n_mn_s) \\ 
\Sigma_{k+1|k+1} = & \quad (\Sigma_{k+1|k}^{-1} + H_k^T\noise^{-1}H_k)^{-1}  & \bigO(n_mn_s^2 ) \label{eqn:filtering}
\end{align}
\end{algorithmic}
\caption{Standard implementation of Kalman Filter}
\label{alg:standard}
\end{algorithm}
Here, in step~\eqref{eqn:residualcovariance}, $S_k$ is known as the innovation covariance and in step~\eqref{eqn:kalmangain}, $K_k$ is known as the Kalman gain matrix. We also define the dimension of the state variables $n_s$ and the number of measurements per time step as $n_m$. We also summarize the computational costs of each step of the Kalman filter, which retains only the leading order terms of $n_s$ with the assumption that the number of measurements per time step $n_m$ is much smaller than the dimension of the state variables $n_s$, i.e. $n_m \ll n_s$. This is typically the case in under-determined inverse problems.

\textbf{Assumptions}: We now state the assumptions that we make in our fast algorithm. The state transition matrix $F_k$ is assumed to be the identity matrix, i.e. $F_k = I$ (see for e.g.~\cite{li2014kalman,soleimani2007dynamic,vauhkonen1998kalman}).  This assumption is known in the literature as the random walk forecast model. This model is useful in practical applications in which data is acquired at a rapid rate when changes in states between times when data is obtained are small that they can approximated by a  random walk process, denoted by $w_k$	. An immediate consequence of this assumption is that it lowers the computational cost in equation for prediction~\eqref{eqn:prediction} from $\bigO(n_s^3)$ to $\bigO(n_s^2)$. Furthermore, we assume that the system noise $w_k$ follows the same distribution $w_k \sim \normal(0,\prior)$ at each step.

Our method will be sufficiently general that we can handle time-varying measurement operator $H_k$. In particular, even the dimensions of $y_k$ will be allowed to vary at each step. However, to describe the general approach it will be convenient to let $H_k=H$. We will relax this assumption later in section. As stated earlier, we also assume that the number of measurements, i.e., the dimension of $y_k$ is small compared to the dimension of the state variable which is assumed to be $\bigO(10^6)$ or higher~\cite{li2014kalman}. 

\subsection{Computational tools}

\subsubsection{Efficient representation of covariance matrices}\label{sec:cov}
The system noise covariance matrix $\prior$ is usually specified in terms of a covariance kernel $\kappa(\cdot,\cdot)$ with entries $\prior(i,j) = \kappa(\bx_i,\bx_j)$. A popular choice for $\kappa(\cdot,\cdot)$ is from the Mat\'{e}rn family of covariance kernels. 
\begin{equation}\label{eqn:maternfamily}
\kappa(\bx,\textbf{y}) = C_{\alpha,\nu}(r) = \frac{1}{2^{\nu-1}\Gamma(\nu)} (\sqrt{2\nu}\alpha r)^\nu K_\nu(\sqrt{2\nu}\alpha r) 
\end{equation}
where $r=\normtwo{\bx-\textbf{y}}$, $\Gamma$ is the Gamma function, $K_\nu(\cdot)$ is the modified Bessel function of the second kind of order $\nu$ and $\alpha$ is a scaling factor. Equation~\eqref{eqn:maternfamily} takes special forms for certain parameters $\nu$. For example, when $\nu=1/2$, $C_{\alpha,\nu}$ corresponds to the exponential covariance function, $\nu = 1/2+n$ where $n$ is an integer, $C_{\alpha,\nu}$ is the product of an exponential covariance and a polynomial of order $n$. In the limit as $\nu\rightarrow\infty$, and for appropriate scaling of $\alpha$, $C_{\alpha,\nu}$ converges to the Gaussian covariance kernel. 

For stationary or translational invariant covariance kernels with points located on a regular equispaced grid, the computational cost for matrix-vector products (henceforth referred to as MVPs) involving the prior covariance matrices can be reduced from $\bigO(n_s^2)$ using the naive approach, to $\bigO (n_s \log n_s)$ by exploiting the connection between Toeplitz structure in 1D or Block-Toeplitz structure in 2D etc, and the Fast Fourier Transform (FFT)~\cite{nowak2003efficient}. For irregular grids, it can be shown that the cost for approximate matrix-vector products (MVPs) involving the prior covariance matrix $\prior$ can be reduced to $\bigO(n_s\log n_s)$ using Hierarchical matrices~\cite{saibaba2012efficient} or $\bigO(n_s)$ using ${\cal{H}}^2$-matrices or the Fast Multipole Method (FMM)~\cite{ambikasaran2012large}. In this work, we will use the Hierarchical matrix approach originally developed by Hackbusch and co-authors~\cite{borm2003introduction,Grasedyck03constructionand} and applied to dense covariance matrices in~\cite{saibaba2012efficient,saibaba2012application,saibaba2013uncertainty,saibaba2013fast}.

For forming MVPs $\prior^{-1}x$, we use an iterative solver such as GMRES with a preconditioner that employs approximate cardinal functions based on local centers and special points~\cite{beatson1999fast}. The cost of constructing the preconditioner is $\bigO(n_s)$ or $\bigO(n_s\log n_s)$. Assuming that the number of iterations is independent of the size of the system, the cost of inverting the prior covariance matrix, i.e. forming $\prior^{-1}x$ is also $\bigO(n_\text{iter}n_s)$ or $\bigO(n_\text{iter}n_s\log n_s)$, where $n_\text{iter}$ is the number of iterations required to converge to the desired tolerance. In conclusion, the cost for forming $\prior x$ and $\prior^{-1} x$ is $\bigO(n_s\log^\gamma n_s)$, where $\gamma \in \{0,1\}$ is a constant depending on the method chosen.

\subsubsection{Low rank representation}\label{sec:lowrank}
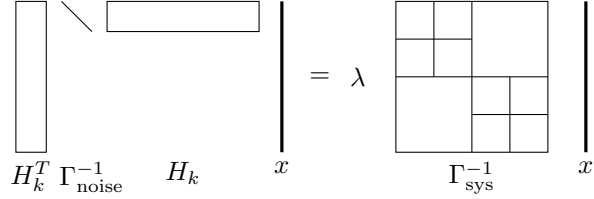
\begin{figure}
\centering
 \begin{tikzpicture}[scale = 0.2]
\pgfmathsetmacro{\H}{25}
\pgfmathsetmacro{\LR}{0}
 
\draw (\LR + 0,0) rectangle (\LR + 2, 10);
\draw (\LR + 3,10) -- (\LR + 5,8); 
\draw (\LR + 6,8) rectangle (\LR + 16,10);
  
\node [below] at (\LR + 1,0) {$H_k^T$};
\node [below] at (\LR + 5,0) {$\noise^{-1}$};
\node [below] at (\LR + 11,0) {$H_k$};
\node (a1) at (\LR +5,10) {};

\draw [very thick] (\LR+17.5,0) -- (\LR + 17.5,10);
\node [below] at (\LR+17.5,0) {$x$};
\node (eq) at (20, 5) {$=$};
\node (lambda) at (22.5,5) {$\lambda$};

\draw (\H,0) rectangle (\H + 10, 10);
\draw (\H,5) -- (\H + 10,5);
\draw (\H+5,0) -- (\H+5,10);
\draw (\H,7.5) -- (\H + 5,7.5);
\draw (\H+2.5,5) -- (\H+2.5, 10);
\draw (\H+5,2.5) -- (\H+10, 2.5);
\draw (\H+7.5,0) -- (\H+7.5,5);
\node (a2) [below] at (\H + 5,0) {$\prior^{-1}$};

\draw [very thick] (\H+12.5,0) -- (\H + 12.5,10);
\node [below] at (\H+12.5,0) {$x$};
\end{tikzpicture}
\caption{Visual representation of the generalized eigenvalue problem described in equation~\eqref{eqn:ghep}. The matrix $\prior$ is assumed to be approximated by a $\mathcal{H}$-matrix.}
\label{fig:ghep}
\end{figure}

We consider the generalized eigenvalue problem (also see figure~\ref{fig:ghep}) 

\begin{equation}\label{eqn:ghep}
H^T\noise^{-1}H x = \prior^{-1} x
\end{equation}

Since both matrices $H_\text{red} \define H^T\noise^{-1}H$ and $\prior$ are symmetric and $ \prior$ is symmetric positive definite, we have the following generalized eigendecomposition
\begin{equation}
\label{eqn:geneigendecomposition}
H_\text{red} = \prior^{-1} U \Lambda U^T \prior^{-1} \qquad U^T\prior^{-1}U = I
\end{equation}

Further properties of this decomposition have also been considered in~\cite{cui2014likelihood,saibaba2013uncertainty}. In general, the rank of this eigendecomposition is $\min\{n_m,n_s\}$. By the assumption that we have made, the number of measurements $n_m$ is much smaller than the number of state variables $n_s$, i.e., $n_m \ll n_s$. This implies that the eigenvalue problem $H_\text{red}x=\lambda \prior^{-1}x$ has a low numerical rank and be efficiently computed. An efficient algorithm for computing the generalized Hermitian eigendecomposition has been proposed in~\cite{saibaba2013randomized}, that avoids forming expensive matrix-vector products with $\prior^{1/2}$ or its inverse $\prior^{-1/2}$ (this is also summarized in the Appendix~\ref{sec:rand}. However, we have chosen to use~\cite[algorithms 3,4]{saibaba2013randomized}. For further details, the reader is referred to~\cite{saibaba2013randomized}. 

For many ill-posed inverse problems, the numerical rank $r$ of the eigenvalue problem, $H_\text{red}x=\lambda \prior^{-1}x$, is small and independent of the problem size, i.e., the number of state variables. The generalized eigendecomposition combines information from the prior and the reduced Hessian and takes advantage of the eigenvalue decay in one (or both) matrices - when the reduced Hessian has rapidly decaying eigenvalues, or the prior is smooth.  Analytical evidence for the eigenvalue decay of the reduced Hessian $H_\text{red}$ is provided in~\cite{flath2011fast} in the context of advection-diffusion based inverse problems and in~\cite{bui2012analysis1,bui2012analysis2} for inverse scattering problems. For the case of system noise covariance matrices $\prior$, the eigenvalue spectrum is known to decay rapidly when the covariance kernels are smooth~\cite{schwab2006karhunen}. The $r$ retained eigenvectors are the modes along which the parameter field is informed by a combination of the data and the prior. Typically the data and prior are informative about the low-frequency modes and as a result local information and fine scale information is hard to recover from the data. In such cases, the rank of the eigendecomposition $r$ satisfies the following inequality $r \leq \min \{n_m,n_s\}$.

\subsection{Fast Kalman Filter for measurement operators $H_k=H$}

In this section, we describe an efficient algorithm for the Kalman filter based on an efficient representation and update of the state covariance matrix $\Skk$. We make the assumption that the measurement operator $H_k$ does not change. In the situation with small number of measurements, we demonstrate that the updates for the Kalman filter can be computed efficiently in $\bigO(n_s)$ or $\bigO(n_s\log n_s)$. Our approach is as follows: we first construct an educated guess for the form that $\Skk$ should take at each time step, and then demonstrate by an inductive argument that the form is preserved at the next time step and can be efficiently calculated. We will consider the constant $H$ case separately since the main idea is easier to explain and this special case is important because the resulting update is numerically exact.

\begin{algorithm}
\begin{algorithmic}[1]
\REQUIRE Measurement operator $H$ and measurements $y_k$ for $k=1,\dots,N_t$, system noise covariance $\prior$ and measurement noise $\noise$ a
\STATE Compute the generalized eigendecomposition \[ H = \prior^{-1}U_k \Lambda_k U_k^T\prior^{-1}\qquad \text{with}\qquad  U_k^T\prior^{-1}U_k = I\]
\STATE Compute $ \Gamma \define \prior H^T$ 
\COMMENT {// Cross Covariance}
\FOR {$k=1,\dots,N_t$}
\STATE Update $\alpha_{k+1}=\alpha_k+1$ and $\tilde{D}_k \define \left((\alpha_{k+1} I - D_k)^{-1} \Lambda^{-1} + I\right)^{-1}$
\STATE Compute $F=\alpha_{k+1}\Gamma - UD_k(U^TH)$
\STATE Compute $\skkpp = \skk + F(HF+\noise)^{-1}(y_k-H\skk)$
\STATE $D_{k+1} \define D_k + \tilde{D}_k^{-1}(\alpha_{k+1}I - D_k) $
\ENDFOR
\end{algorithmic}
\caption{Fast Kalman Filter for random walk forecast model}
\label{alg:kalman}
\end{algorithm}

We start with the following ansatz for an efficient representation of the a posteriori estimate covariance matrix $\Sigma_{k|k}$
\begin{equation}\label{eqn:sigmakk}
\Sigma_{k|k} = \alpha_k\prior - W_kD_kW_k^T  
\end{equation}
and the matrices $W_k$ are chosen to be the eigenmodes of the generalized eigenvalue problem described in equation~\eqref{eqn:ghep}, i.e. $W_k = U$. This representation (also see figure~\ref{fig:lowrankperturb}) assumes that the a posteriori estimate covariance matrix $\Sigma_{k|k}$ can be written as a weighted combination of two terms - the system noise covariance matrix and a low-rank term term that contains eigenvectors of the generalized eigenvectors of the eigenvalue problem described in~\eqref{eqn:geneigendecomposition}. We will show that with this ansatz, the a priori estimate covariance matrix $\Sigma_{k+1|k}$ and a posteriori estimate covariance matrix $\Sigma_{k+1|k+1}$ can be written in a similar form  as equation~\eqref{eqn:sigmakk} with updated weights $\alpha_{k+1}$ and $D_{k+1}$. The advantage with this is representation is twofold, (1) it provides an efficient representation of $\Skk$ since $\prior$ can be efficiently using, for example, the ${\cal H}$-matrix approach in $\nlogn$ or $\bigO(n_s)$\cite{saibaba2012efficient,saibaba2012application,ambikasaran2012large} and (2) as we will show, the weights $\alpha_{k+1}$ and $D_{k+1}$ can be efficiently updated in a cost $\bigO(r)$, where $r$ is the rank of the low-rank representation. 

\begin{figure} \centering
  \begin{tikzpicture}[scale = 0.15]
\pgfmathsetmacro{\H}{17}
\pgfmathsetmacro{\LR}{32}
   \draw (0,0) rectangle (10,10);
    \node (A) at (5,5) {$\Skk$};
 \node (eq) at (12, 5) {$=$};

\node (alpha) at (15,5) {$\alpha_k$};

\draw (\H,0) rectangle (\H + 10, 10);
\draw (\H,5) -- (\H + 10,5);
\draw (\H+5,0) -- (\H+5,10);
\draw (\H,7.5) -- (\H + 5,7.5);
\draw (\H+2.5,5) -- (\H+2.5, 10);
\draw (\H+5,2.5) -- (\H+10, 2.5);
\draw (\H+7.5,0) -- (\H+7.5,5);
\node (p) at (\H + 13,5) {$-$};

\node [below] at (\H + 5,0) {$\prior$};
\draw (\LR + 0,0) rectangle (\LR + 2, 10);
\draw (\LR + 3,10) -- (\LR + 5,8); 
\draw (\LR + 6,8) rectangle (\LR + 16,10);
  
\node [below] at (\LR + 1,0) {$W_k$};
\node [below] at (\LR + 4,8) {$D_k$};
\node [below] at (\LR + 11,8) {$W_k^*$};
\end{tikzpicture}
\caption{Visual representation of the efficient representation of the state covariance matrix $\Skk$ as the weighted sum of the system noise covariance matrix $\prior$ and a low-rank perturbation term}
\label{fig:lowrankperturb}
\end{figure}

Plugging the ansatz of equation~\eqref{eqn:sigmakk}, into the prediction equation~\eqref{eqn:prediction}, we obtain the following representation
\begin{align*}
\Sigma_{k+1|k} \quad = & \quad   \alpha_k\prior - UD_kU^T  + \prior \\ 
 =& \quad   \alpha_{k+1}\prior - UD_kU^T    
\end{align*}
Recall that the above relation holds because $F_k = I$, which corresponds to the random walk forecast model. Further, observe that the covariance matrix $\Sigma_{k+1|k}$ is of the same form as our assumption for $\Sigma_{k|k}$. Now, consider the filtering equation in~\eqref{eqn:prediction}. Plugging in the representation for $\Sigma_{k+1|k}$, we have 
\[ \Sigma_{k+1|k+1} =  (\Sigma_{k+1|k}^{-1} + H^T\noise^{-1}H)^{-1} = \quad ( I + \Sigma_{k+1|k}  H^T\noise^{-1}H)^{-1} \Sigma_{k+1|k} \]
Next, we plug in the eigendecomposition in~\eqref{eqn:geneigendecomposition} into the above expression to yield
\begin{align*}
\Sigma_{k+1|k}  H^T\noise^{-1}H  \quad =& \quad   \Sigma_{k+1|k}   \prior^{-1} U \Lambda U^T \prior^{-1}  \\ 
= & \quad \alpha_{k+1} U\Lambda U^T\prior^{-1} - U\Lambda D_k U^T\prior^{-1} \\
= & \quad U\Lambda (\alpha_{k+1} I - D_k)U^T\prior^{-1}  
\end{align*}
We have used the identity that $U^T\prior^{-1}U = I$ by construction. Applying the Sherman-Morrison-Woodbury identity to $( I + \Sigma_{k+1|k}  H^T\noise^{-1}H)^{-1} $ gives us 
\begin{align*}
( I + \Sigma_{k+1|k}  H^T\noise^{-1}H)^{-1}  \quad = & \quad  (I + \quad U\Lambda (\alpha_{k+1} I - D_k)U^T\prior^{-1}  )^{-1} \\ 
= & \quad I -  U( (\alpha_{k+1} I - D_k)^{-1} \Lambda^{-1} + I)^{-1}U^T \prior^{-1} \\ 
 = & \quad I - U\tilde{D}_k^{-1}U^T\prior^{-1}  
\end{align*}
where $\tilde{D}_k^{-1} \define (\alpha_{k+1} I - D_k)^{-1} \Lambda^{-1} + I$. Finally, 
\begin{align*}
\Sigma_{k+1|k+1} \quad = & \quad ( I + \Sigma_{k+1|k}  H^T\noise^{-1}H)^{-1} \Sigma_{k+1|k}  \\ \nonumber
= &\quad \left(I -  U\tilde{D}_k^{-1}U^T \prior^{-1} \right) (\alpha_{k+1}\prior - UD_kU^T) \\  
= & \quad  \alpha_{k+1} \prior - \alpha_{k+1}U\tilde{D}_k^{-1}U - UD_kU^T + U\tilde{D}_k^{-1}D_kU^T \\ 
= & \quad \alpha_{k+1} \prior - U D_{k+1} U^T  
\end{align*}
where $D_{k+1} \define D_k + \tilde{D}_k^{-1}(\alpha_{k+1}I - D_k) $ is the updated coefficient. Note that despite the complicated expression for $D_{k+1}$, it is still a diagonal matrix. The inductive proof is completed by making the assumption that $\Sigma_{0|0}= \alpha_0\prior$ which is trivially of the form in equation~\eqref{eqn:sigmakk}.

The algorithm for computing the updates $\skk$ and $\Skk$ is summarized in algorithm~\ref{alg:kalman} and the relevant computational costs at each time step are discussed in algorithm~\ref{alg:kalmanconstantH}.

\begin{algorithm}
\begin{algorithmic}
\STATE\textbf{Predict}
\begin{align}
\hat{s}_{k+1|k} \quad = & \quad \hat{s}_{k|k} & - \\  
\alpha_{k+1} \quad = & \quad \alpha_{k} + 1 & \bigO(1)  
\end{align}

\STATE \textbf{Update}
\begin{align}
\Skkp H^T \quad  = & \quad \alpha_{k+1}\prior H^T - UD_k(HU)^T & \bigO(rn_mn_s) \\
S_k \quad = &  \quad   \alpha_{k+1}H\prior H^T - HUD_k(HU)^T  + \noise & \bigO(rn_m)  \\ 
\skkpp \quad = & \quad \skkp+ \Skkp H^T S_k^{-1}( y_k - H\skk ) &  \bigO(n_mn_s) \\ 
D_{k+1}\quad  = & \quad  D_k + \left(  (\alpha_{k+1} I - D_k)^{-1} \Lambda^{-1} + I \right)^{-1}(\alpha_{k+1}I - D_k) & \bigO(r) \\
\end{align}
\end{algorithmic}
\caption{Fast Kalman Filter for Random Walk Forecast Model}
\label{alg:kalmanconstantH}
\end{algorithm}

\subsection{The case where $H_k \neq H$}\label{sec:nonconstantH}

In this section, we will relax the assumption that the measurement operator is constant at each step. Essentially, we now need to compute a new low-rank decomposition $H_k^T\noise^{-1}H_k = \prior^{-1}U_k\Lambda_kU_k^T\prior^{-1}$ at each iteration. This can be efficiently computed using the randomized algorithm for Generalized Hermitian Eigenvalue problem described in section~\ref{sec:lowrank} and in the reference~\cite{saibaba2013randomized}. The measurement operator changes at each iteration; consequently, the eigenmodes along which information propagates change from iteration to iteration. As a result we expect the rank of the low-rank perturbation to grow.  

Let us begin by making the ansatz that $\Skk = \alpha_k\prior - W_kD_kW_k^T$ where $W_k^T\prior^{-1}W_k = I$ and as before $\alpha_k$ is a scalar and $D_k$ is a diagonal matrix. $W_k$ are no longer simply the generalized eigenvectors but are now updated at each step. Since the measurements do not enter into the prediction, the expression for $\Skkp$ is similar to equation~\eqref{eqn:sigmakk} and is given by $\Skkp=\alpha_{k+1}\prior - W_kD_kW_k^T$. Now we use the Sherman-Morrison-Woodbury update to derive an expression for $\Skkp^{-1}$  

\begin{align*}
\Skkp^{-1} \quad = & \quad  \alpha_{k+1}^{-1} \left(  \prior - \alpha_{k+1}^{-1}W_kD_kW_k^T \right) \nonumber \\ 
 = &\quad \alpha_{k+1}^{-1} \left( \prior^{-1} - \prior^{-1}W_k (I - \alpha_{k+1}D_k^{-1})^{-1}W_k\prior^{-1}\right) \nonumber \\ 
= & \quad \alpha_{k+1}^{-1}  \prior^{-1} - \prior^{-1}W_k \alpha_{k+1}^{-1}( \alpha_{k+1}I + D_k)^{-1}D_kW_k\prior^{-1} \nonumber
\end{align*}

We seek an efficient representation for the matrix $\Skkpp$, which we write as 
\begin{align*}
\Skkpp^{-1} \quad = & \quad \Skkp^{-1} + H_k^T\noise^{-1} H_k \nonumber\\ 
= & \quad \alpha_{k+1}^{-1}  \prior^{-1} - \prior^{-1}W_k \alpha_{k+1}^{-1}( \alpha_{k+1}I - D_k)^{-1}D_kW_k\prior^{-1}  + \prior^{-1} U_k\Lambda_k U_k^T \prior^{-1}\nonumber\\ 
= & \quad  \alpha_{k+1}^{-1}  \prior^{-1} + \prior^{-1} \left(-W_k \alpha_{k+1}^{-1}( \alpha_{k+1}I - D_k)^{-1}D_kW_k +  U_k\Lambda_k U_k^T  \right) \prior^{-1}
\end{align*}

The last line in the above equation is the sum of two low-rank matrices. In general, with the addition of two low-rank matrices, the rank of their sum is bounded by the sum of their ranks. In order to make the low-rank terms orthogonal with respect to the matrix $\prior^{-1}$, we apply algorithm~\ref{alg:lowrankaddition} (described in the Appendix~\ref{sec:addlowrank}) as 
\[ [\hat{W},\hat{D}] = \text{AddLowRank}\left(W_k,\bar{D}_k, U_k,\Lambda_k,\prior^{-1}\right)\]
where for convenience we define $\bar{D}_k \define -\alpha_{k+1}( \alpha_{k+1}I - D_k)^{-1}D_k$. Finally, we can express apply the Sherman-Morrison-Woodbury identity to the matrix $\Skkpp$ as 
\begin{align*}
\Skkpp^{-1} \quad = & \quad \alpha_{k+1}^{-1}  \prior^{-1} + \prior^{-1}\hat{W}_k\hat{D}_k\hat{W}_k^T \prior^{-1} \nonumber \\
\Skkpp \quad  = & \quad \alpha_{k+1} \prior - \hat{W}_k\alpha_{k+1} (I + \alpha_{k+1}^{-1}\hat{D}_k^{-1})^{-1} \hat{W}_k^T \nonumber \\
= & \quad \alpha_{k+1} \prior - W_{k+1} D_{k+1} W_{k+1}^T \nonumber
\end{align*}
where $W_{k+1} \define \hat{W}_k$ and $D_{k+1} \define \alpha_{k+1} (I + \alpha_{k+1}^{-1}\hat{D}_k^{-1})^{-1} $. We have therefore shown that every subsequent iterate using the Kalman update can be written in the form of the ansatz $\Skk = \alpha_k\prior - W_kD_kW_k^T$, where $D_k$ is a diagonal matrix and $W_k^T\prior^{-1}W_k = I$. However, observe that the rank $r_{k+1}$ of the low-rank part $\Skkpp$ is now at most $r_{k+1} = r_k + n_{m,k}$. Based on this update, it is easy to see that $r_k = \sum_{l=1}^{k}n_{m,l}$. The update for $\Skk$ is efficient as long the rank $r_k \ll n_s$. One possibility is that the number of measurements per time step is small as well as the number of time steps over which the Kalman filter is applied. Another practical situation is that the spectrum of the GHEP defined in equation~\eqref{eqn:ghep} (that combines the spectral decay in the measurement operator and system noise covariance matrix)  decays rapidly. Furthermore, with each iteration the eigenmodes of the GHEP may be aligned; as a result, the effective rank of the low-rank perturbation grows at a slower rate than the number of measurements up to a given time. This is resulted by the results in section~\ref{sec:extkalman}. In summary, the computational cost of updating a step is $\bigO(n_s\log n_s + r^2n_s)$, where $r$ is the rank of the perturbation.


\subsection{Large-scale implementation}\label{sec:parallel}
 Although our algorithm has an asymptotic computational complexity of $\bigO(N\log N)$, scalable implementations are crucial for large-scale systems arising from finely discretized grids that are needed for accurate simulations of real world applications.   In Algorithm~\ref{alg:kalman} the most  expensive steps are those that involve the system noise covariance matrix $\prior$. Specifically, these involve  the computation of the generalized eigenmodes $U$ and the cross-covariance matrix $\prior H^T$. In~\cite[Chapter 2.5]{saibaba2013fast} we demonstrated scalability of matvecs involving $\prior$ systems of sizes which are $\bigO(10^5 )$. For larger problems, there is a need to turn to distributed computing to handle this computational burden. 
 
 There are several parallel algorithms  for efficiently representing $\mathcal{H}$-matrices and computing matvecs of the form $\prior x$. We mention a few relevant references here which are not exhaustive.  
 A highly scalable implementation of the algorithms for $\mathcal{H}$- and $\mathcal{H}^2$ matrices with computational cost for matvecs $\bigO\left(\frac{rN\log N}{p} + r\log^2p\right)$, $ \frac{rN\log N}{p}$ memory usage and $r^2 \log p$ communication costs is publicly available in~\cite{poulson2013dmhm}. Here $r$ is the blockwise rank and $p$ are the number of processors.  This code has been shown to scale up to $\bigO(10^3)$ processors. A related class of matrices called the Hierarchical Semiseparable (HSS) matrices has been developed and algorithms for matvecs of the form $\prior x$, in which $\prior$ is represented as an HSS matrix, have been shown to scale up to $6.4$ billion unknowns on $4096$ processors~\cite{wang2012efficient}. The Fast multipole method is another related algorithm for efficiently computing matvecs of the form $\prior x$, and has also been shown to scale on heterogeneous architectures (CPUs and GPUs) up to 64k cores and 30 billion unknowns~\cite{lashuk2012massively}. The relation between the $\mathcal{H}$-matrices, HSS approach and FMM has been previously discussed in~\cite{ambikasaran2013fast}. Any one of these implementations can be used to handle covariance matrix calculations on a massive scale. 
 
 Two other modifications make our algorithm more efficient in distributed computing. First, the computation of the cross-covariance matrix $\prior H^T$ can be trivially parallelized since the matvec computations are independent of one another.  Second, the randomized algorithm that we are using to compute the dominant eigenmodes is advantageous over Krylov subspace methods for distributed computing since it allows to organize our computations to exploit parallelism. This has been discussed in greater detail in~\cite{saibaba2013randomized}. To summarize, the use of distributed computing in combination with $\mathcal{H}$-matrices can be used for large-scale implementations of Kalman filters. However, in this paper, we restrict our discussion to the algorithm for a single processor. 

\section{Uncertainty quantification}\label{sec:uncert}
We have demonstrated that a low-rank perturbative approach for the state covariance matrix $\Skk$ leads to an efficient representation as well as an efficient algorithm for updating the state estimate. We further demonstrate that this efficient representation is further useful in computing measures of uncertainty of the distribution. In this section, we will assume that $\Skk = \alpha_k\prior - W_kD_kW_k^T$ where $W_k^T\prior^{-1}W_k = I$. Of course, when the observation operator is time invariant, i.e., $H_k=H$, we have that the set of vectors $W_k = U$, where $U$ are the eigenvectors of the eigenvalue problem $H^T\noise^{-1}Hx = \lambda \prior^{-1}x$.  

\subsection{Variance computation}

The variance of the distribution $s_k \sim {\cal N} (\skk, \Skk)$ is obtained by computing the diagonals of the matrix $\Skk$. This can be efficiently done by considering the representation of $\Skk$ by the formula~\eqref{eqn:sigmakk}. Therefore, 
\[\text{Var} (\Skk) \quad \define \quad  \alpha_k \text{diag} (\prior) - \text{diag} (W_kD_kW_k^T)  \]
The resulting computations for the variance scales as $\bigO(rn_s)$. A related measure of uncertainty is the A-optimality criterion, which takes the form $\phi_A \define \frac{1}{n_s}\text{Trace}(A\Skk)$. For $A=I$, the expression for the A-optimality criterion simplifies to \[\phi_I = \trace\Skk = \alpha_k \trace( \prior) - \trace (W_kD_kW_k^T)\]
For $A\neq I$, the A-optimality criterion can be approximated using the Hutchinson trace estimate. For further details, the reader is referred to~\cite{saibaba2013uncertainty}. 

\subsection{Entropy and Relative entropy}
The entropy of a random variable $X$ with probability density function $p(X)$ is defined as $H[X] =  E[ -p(X)\log p(X)] $, where $E[\cdot]$ is the expectation. For Gaussian distributions $X \sim \normal(\mu,\Sigma)$, the entropy can be calculated as $H[X] = \frac{1}{2}\log 2\pi e + \frac{1}{2} \log \det (\Sigma)$. For the distribution $s_k \sim {\cal N} (s_{k|k}, \Skk)$, we entropy can be calculated as 

\[ H[s_k]  = \frac{1}{2}\log 2\pi e + \frac{1}{2} \logdet \Skk  \]
Now consider $\logdet \Skk$. By factorizing out $\prior$ we have 
\begin{align*}
\logdet \Skk \quad  = & \quad \logdet (\alpha_k \prior - UD_kU^T) \\ 
= & \quad \logdet \prior + \logdet (\alpha_k I - \prior^{-1}UD_k U^T) \\ 
= & \quad \logdet \prior + \logdet (\alpha_k I - U^T\prior^{-1}U D_k) \\ 
= & \quad \logdet \prior +  \logdet (\alpha_k I -  D_k) 
\end{align*}
where in the penultimate step we have use Sylvester's determinant lemma and used the orthogonality of $U$ from equation~\eqref{eqn:geneigendecomposition}. Computing $\logdet \prior$ can be computationally challenging but observe that the entropy of $w_k$ is $H[w_k]=   \frac{1}{2}\log 2\pi e + \frac{1}{2} \logdet \prior$. Instead we compute relative entropy, defined between random variables $X$ and $Y$ as $H[Y|X] = H[Y] - H[X]$. It can be readily seen that 
\begin{equation}
H[s_{k|k}|w_k ]\quad  =  \quad \frac{1}{2}\logdet (\alpha_k I -  D_k) 
\end{equation}

\subsection{Sampling from the distribution $s_{k} \sim {\cal N} (\skk, \Sigma_{k|k})$}\label{sec:sampling}
By the modeling assumptions that we have made thus far, the system state vector is normally distributed at each iteration $k$ with the distribution completely specified by the mean $\skk$ and covariance $\Skk$. However, it is difficult to visualize this large matrix. We have shown previously that the statistics of this distribution can be summarized using a few notions of uncertainty. An alternative way to make sense of this distribution, to study and visualize the uncertainty, is to draw samples from the distribution corresponding to the current estimate of the state, namely $\normal (\skk, \Sigma_{k|k})$. Drawing conditional samples is extremely beneficial in analyzing different possible scenarios which are equally ``likely''.  

In general, in the absence of an efficient representation for $\Skk$, sampling from this distribution requires the computation of the Cholesky factorization of a dense matrix. This is computationally expensive as the computational cost scales as $\bigO(n_s^3)$. However, in this section we propose an efficient sampling technique based on the efficient representation of the state covariance matrix $\Skk$. We follow the approach described in~\cite{bui2013computational}. 

Suppose there exists a factorization of $\Sigma_{k|k} = LL^T$, then the samples can be computed as 
\[ s_{k|k} = \hat{s}_{k|k} + Ls_u \qquad s_u \sim {\cal N}(0, I)\] 
It remains to be shown how to construct such a factorization. We can factorize  $\Sigma_{k|k}$  as 
\[ \alpha_k^{-1}\Sigma_{k|k} =  \prior - \alpha_k^{-1}UD_k U^T = \prior^{1/2} \left( I - \alpha_{k+1}^{-1}\prior^{-1/2}UD_kU^T  \prior^{-1/2}\right) \prior^{1/2}\] 
Defining $W \define \prior^{-1/2}U_k $ and with the observation that the columns of $W$ are orthonormal, we consider the matrix $I - WD_k W^T $ which is a low-rank perturbation of the identity matrix
\[ I - WD_k W^T = (I - W\Sigma W^T) (I - W\Sigma W^T) \]
A straightforward calculation shows that the above identity holds true when $\Sigma$ satisfies the matrix quadratic equation
\[\Sigma^2 - 2\Sigma + \alpha_{k}^{-1}D_k = 0 \qquad\Rightarrow \qquad  \Sigma_{\pm} = I \pm (I - \alpha_{k}^{-1}D_k)^{1/2} \] 
With $\Sigma$ computed as above, define $L$ as 
 \begin{align*}
 L \quad \define & \quad \alpha_{k}^{1/2}\prior^{1/2} (I - W\Sigma_\pm W^T) \\
 = & \quad  \alpha_{k}^{1/2}\prior^{1/2} - \alpha_{k}^{1/2}U\Sigma_\pm U^T\prior^{-1/2}
\end{align*}
It can be readily verified that $L$ satisfies $ \Sigma_{k|k}  = LL^T$, and is therefore an approximate square-root of $\Sigma_{k|k}$.

As can be readily seen, sampling from $\normal(\skk, \Sigma_{k|k})$ requires products with $\prior^{1/2}$ and $\prior^{-1/2}$. In other words, it requires unconditional samples from the distributions $\normal(0,\prior)$ and $\normal(0,\prior^{-1})$. Several matrix-free techniques exist in the literature for computing matrix-vector products (MVPs) $\prior^{1/2}x$ and $\prior^{-1/2}x$, that are based on polynomial approximation~\cite{chen2011computing,dietrich1995efficient} or rational approximations and contour integrals~\cite{hale2008computing}. However, the convergence of polynomial approximations is only algebraic when the smallest eigenvalue is close to zero. Rational approximations and contour integral based methods  do not suffer from the same problem, however they require solutions of a number of shifted systems totaling $\bigO (\log\kappa(\prior)) $, where $\kappa(\cdot)$ is the condition number. Although the number of systems to be solved is often small, even for ill-conditioned problems, solving each system can be expensive in practice.

If we restrict ourselves to understanding how uncertainty propagates by studying how the conditional sample corresponding to the same realization change over time, then we need to compute MVPs  $\prior^{1/2}s_u$ and $\prior^{-1/2}s_u$ only once per realization. As a result, it can be treated as a pre-computation. The ability to efficiently propagate realizations makes it similar to the Ensemble Kalman Filter, in which an ensemble of realizations are propagated at each iteration. In this sense, we claim that our approach combines the optimality (in terms of accuracy) of the original Kalman Filter with the computational advantages of the ensemble based approach.


\section{Numerical Experiments}\label{sec:numerical}
The algorithms were implemented in Python using NumPy~\cite{van2011numpy} and SciPy~\cite{jones2001scipy} packages. All the figures were generated using Matplotlib~\cite{Hunter:2007}. 
\subsection{Application: CO$_2$ monitoring}
\subsubsection{Problem setup}\label{sec:setup}
In this application, we consider a synthetic setup of cross-well tomography. The goal is to the image the slowness in the medium where slowness is defined as the reciprocal of seismic velocity. A detailed reservoir model was built for the Frio-II brine pilot CO$_2$ injection experiment using TOUGH2~\cite{pruess1999tough2}. CO$_2$ was injected into a brine aquifer and the simulation predicted the spatial distributions of CO$_2$ and pressure over $5$ days. CO$_2$ can be monitored seismically by mapping the time-varying CO$_2$-induced velocity reductions from measurements of travel-time delays. Each source-receiver pair generates one measurement and therefore, there are $n_y = n_\text{rec} n_\text{sou}$ measurements. Here $n_\text{rec}$ is the number of receivers and $n_\text{sou}$ are the number of sources. In this application, we pick $n_\text{sou} = 6$ and $n_\text{rec} = 48$ and synthetic travel-time delay measurements are obtained every 3 hours. The domain is discretized into a sequence of grids of sizes $59\times 55, 117\times 109$ and $234\times 219$. Within each cell, the slowness is assumed to be constant. Therefore, the time taken from the source to the receiver is a weighted sum of the slowness in the cell, weighted by the length of the ray within the cell. The acquisition geometry is the same as in~\cite{li2014kalman} and remains fixed during the monitoring experiment.

The travel times are obtained by integrating the slowness along the ray path on which the seismic waves propagate. As a first order approximation, the seismic wave is modeled as traveling along a straight line from the sources to the receivers without reflections or refractions. the measurement takes the following form 

\begin{equation} \label{eqn:lineartomography}
y_t = \int_l s(r) dl \approx Hs_t 
\end{equation}
where $y_t$ are the observed (synthetic) travel times, $s$ is the slowness that we are interested in imaging and $H$ is the measurement operator, whose rows correspond to each source-receiver pair and are constructed such that their inner product with the slowness would result in the travel time. It is represented as a sparse matrix with $\bigO (n_m\sqrt{n_s})$ non-zero entries - each row has $\bigO (\sqrt{n_s})$ entries and there are $n_m$ measurements. CO$_2$-induced low velocity zone is imaged from travel-time delays $\Delta y_t$ relative to the baseline travel-time, which can be obtained by subtracting the baseline travel-time from the measurement equation~\ref{eqn:lineartomography}. That is, $\Delta y_t = H\Delta s_t$. The variable of interest $\Delta s_t$ is the perturbation of the background slowness at time step $t$. The differential tomography approach applies spatial and temporal regularizations directly on the slowness perturbations instead of slowness itself. Further details about the synthetic setup can be found in the following references~\cite{ambikasaran2012large,li2014kalman} and will not be described here.

\begin{figure} \centering
\begin{tikzpicture}[scale = 0.3]
\draw (0,0) rectangle (10,10);
\foreach \y in {2.5,5.5} {
	\node at (0.5, \y) [circle, fill = red, scale = 0.4] {};
}
\foreach \y in {2.5,4.5,6.5,8.5} {
	\node at (9.5, \y) [circle, fill = blue, scale = 0.4] {};
}

\draw (0.5,5.5) -- (9.5,2.5);
\draw (0.5,2.5) -- (9.5,8.5);

    \foreach \i in {\xMin,...,\xMax} {
        \draw [very thin,gray] (\i,\yMin) -- (\i,\yMax)  node [below] at (\i,\yMin) {};
    }
    \foreach \i in {\yMin,...,\yMax} {
        \draw [very thin,gray] (\xMin,\i) -- (\xMax,\i) node [left] at (\xMin,\i) {};
    };

\node (s) at (12,7) [circle, fill = red, scale = 0.4] {};
\node (r) at (12,3) [circle, fill = blue, scale = 0.4] {};
\node at (15, 7)[right of = s] {Sources};
\node at (15, 3) [right of = r] {Receivers};
\end{tikzpicture}
\caption{Visual representation of the tomographic setup for monitoring CO$_2$ concentration described in section~\ref{sec:setup}.}
\label{fig:tomosetup}
\end{figure}
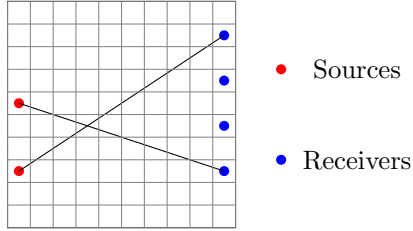

Other parameters chosen for the reconstruction are as follows. The covariance kernel is chosen to be $\kappa(r) = \theta\exp\left(-\frac{r^p}{l^p} \right)$, with $p = 1/2$ and $\theta = 10^{-4}$ . We also assume that $\noise = \sigma^2I$ with $\sigma^2 = 2\times 10^{-4}$. The same noise is added to the measurements. We assume that no CO$_2$ is present before the injection, and as a result we assume that $s_0=\Sigma_{0|0} = 0$.

\subsubsection{Assimilation results}

\begin{figure}[!ht]
\centering
\includegraphics[scale=0.28]{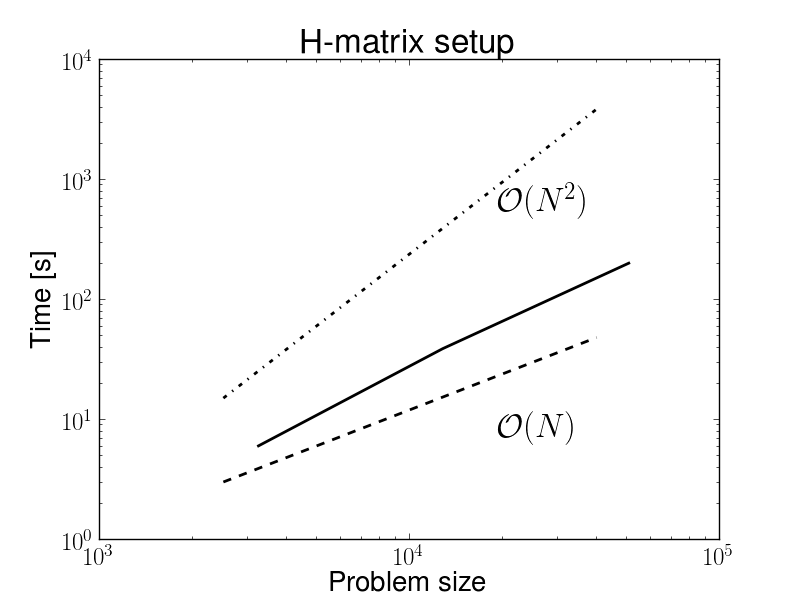}
\includegraphics[scale=0.28]{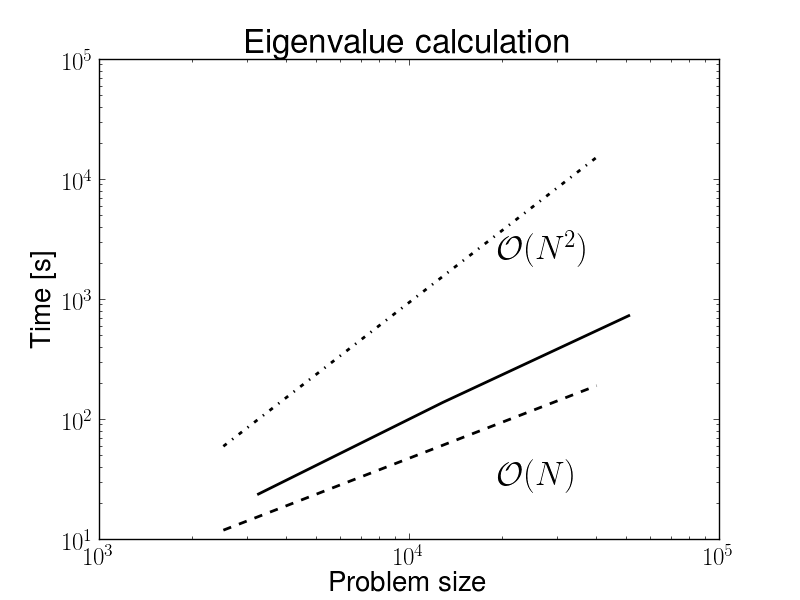}
\caption{(left) Cost for constructing the ${\cal{H}}$-matrix and (right) cost of computing the dominant eigenmodes of the generalized eigenvalue problem $H^T\noise^{-1}H x = \lambda \prior^{-1} x$. The number of measurements are $288$ and the grid sizes varied from $59\times 55$ to $234\times 219$. A block tolerance of $10^{-8}$ was used in the construction of the $\mathcal{H}$-matrix-matrix. See~\cite{saibaba2012application} for further details.  } 
\label{fig:offlinecosts}
\end{figure}

We now discuss the  computational costs associated with performing the assimilation for each measurement time step for various problem sizes for grids varying from  $59\times 55, 117\times 109$ and $234\times 219$. The Fast Kalman Filter that we propose has computational costs which include an offline stage and an online stage. In the offline stage, the computational costs consist of the time spent in constructing the Hierarchical matrix and the time spent in computing the eigenmodes of the generalized eigenvalue problem $H^T\noise^{-1}H x = \lambda \prior^{-1} x$. The cost for constructing the ${\cal{H}}$-matrix scales as $\bigO(k^2n_s)$ where $k$ is the block-wise rank chosen such that the relative Frobenius norm is $\varepsilon$. Further, since computing matrix-vector products with $\prior$ and $\prior^{-1}$ scale as $\bigO(kn_s)$ and $H$ is a sparse matrix, the resulting cost of computing the dominant eigenmodes also scale as $\bigO(n_s)$. This is confirmed by the plots in figure~\ref{fig:offlinecosts}. The online computational costs of the fast Kalman filter is also $\bigO(n_s)$ and is summarized in algorithm~\ref{alg:kalmanconstantH}. 

\begin{figure}[!ht]
\centering
\includegraphics[scale=0.5]{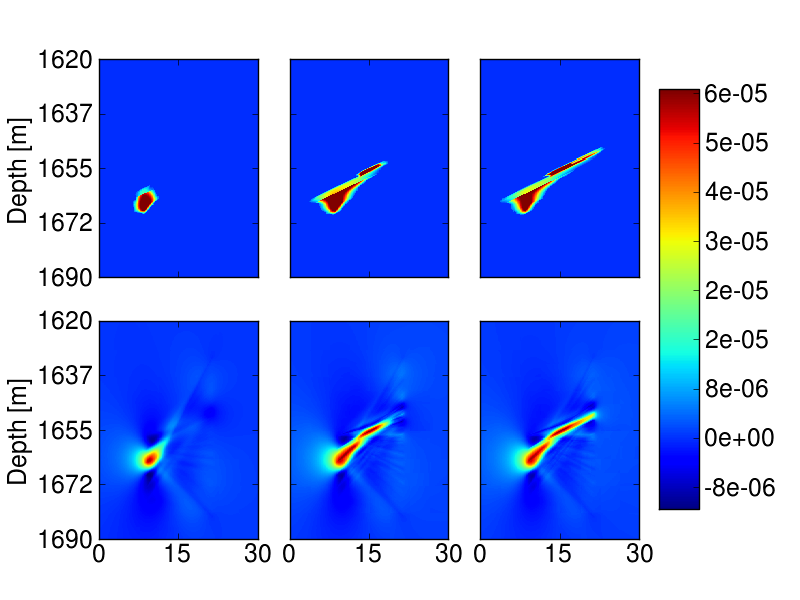}
\caption{True (above) and estimated (below) CO$_2$-induced changes in slowness (reciprocal of velocity) between two wells for the grid size $234\times 219$ (finest grid) at times $3$, $30$ and $60$ hours respectively.}
\end{figure}

\begin{figure}[!ht]
\centering
\includegraphics[scale=0.29]{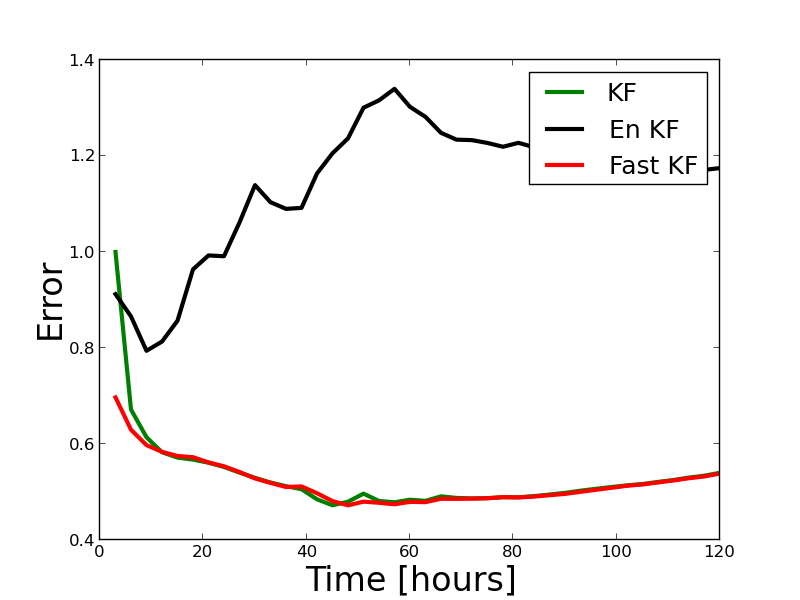}
\includegraphics[scale=0.29]{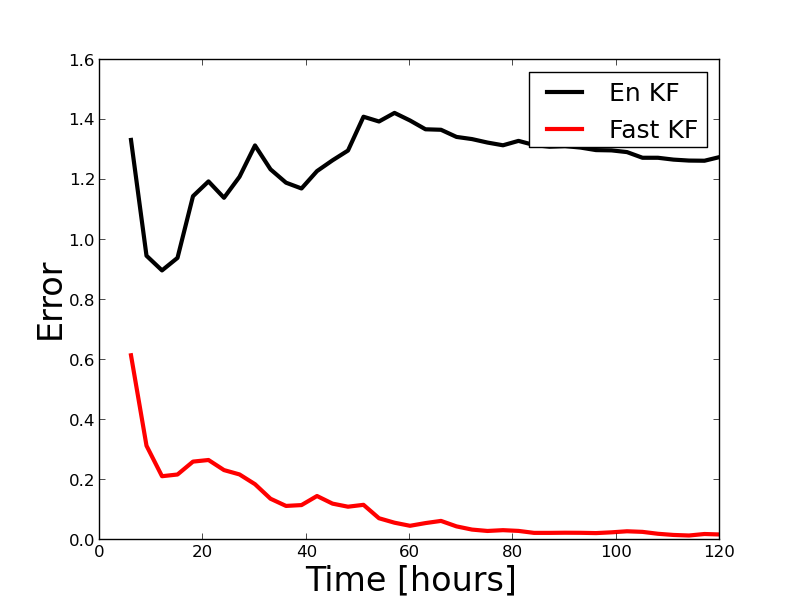}
\caption{(left) Error of the computed solution using the fast algorithm compared against the solution from EnKF and the exact Kalman Filter. The parameters of the covariance kernels and noise are defined in section~\ref{sec:setup}.  For the Ensemble Kalman Filter, $1000$ ensembles were used. (right) Errors in the reconstruction of the proposed fast algorithm and the Ensemble Kalman filter. Instead of the true solution we use the reconstruction from the standard Kalman filter (with the same parameters) as the true solution. }
\label{fig:errkalman}
\end{figure}

\begin{figure}[!ht]
\centering
\includegraphics[scale=0.28]{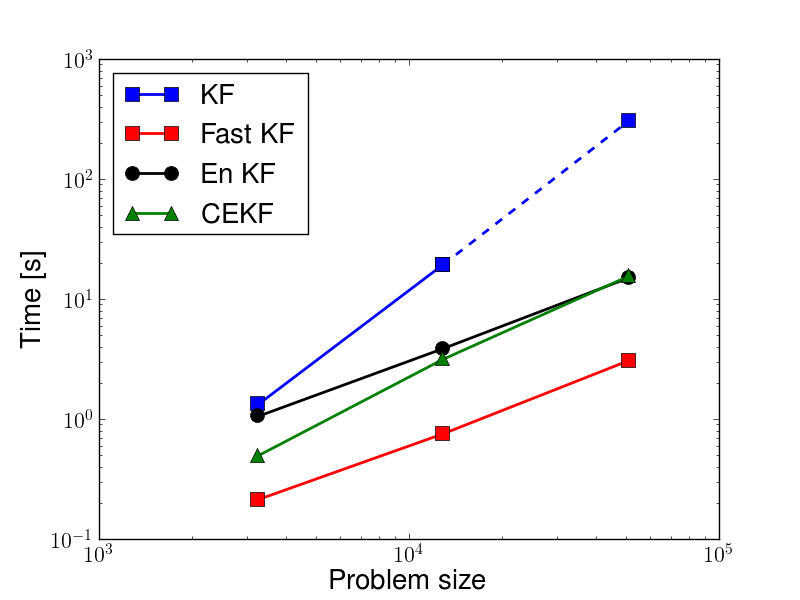}
\includegraphics[scale=0.28]{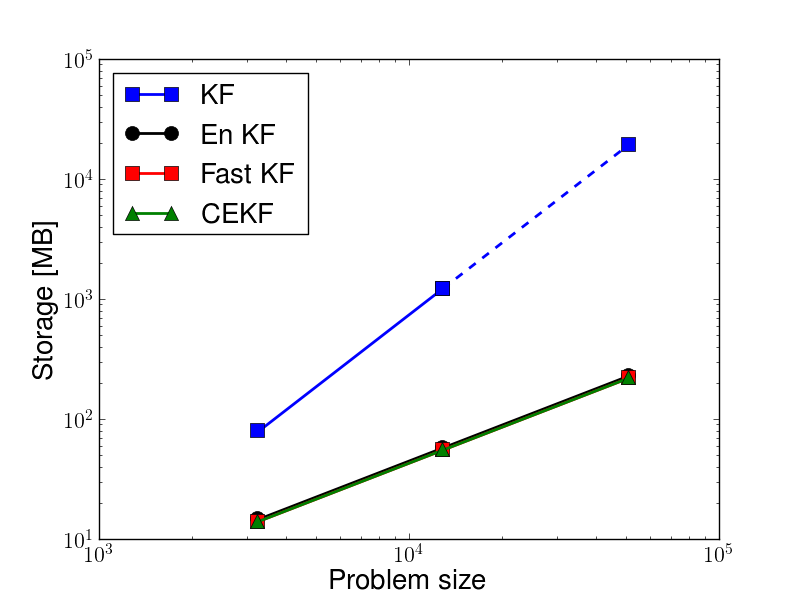}
\caption{The computational cost (left) and storage cost (right) for assimilating measurements for a single assimilation step for  the standard Kalman filter (KF), Ensemble Kalman Filter (EnKF) and the fast Kalman Filter for random walk proposed in~\cite{li2014kalman} (CEKF) and the fast Kalman Filter proposed in this paper (FKF). The full Kalman algorithm was not run on the grid $234\times 219$ and the dashed line indicates the expected time calculated by extrapolation.}
\label{fig:costkalman}
\end{figure}

Comparisons are performed against the standard Kalman filter (KF), Ensemble Kalman Filter (EnKF)\footnote{The details of the implementation of the Ensemble Kalman Filter that we use in this paper are provided in~\cite{li2014kalman}} and the fast Kalman Filter proposed in~\cite{li2014kalman} (CEKF) and the fast Kalman Filter proposed in this paper (FKF). The errors are displayed in figure~\ref{fig:errkalman} and the computational costs are displayed in figure~\ref{fig:costkalman}. As mentioned earlier, the storage and computational costs of the standard implementation of the Kalman Filter scales as $\bigO(n_s^2)$. At the finest scale, the cost for a single assimilation step using the standard Kalman Filter requires over $4$ hours, whereas using our fast algorithm it requires only a few seconds. The computational and storage cost of EnKF also scale linearly with the number of unknowns, as it propagates errors using an ensemble consists of $N$ realizations of state vectors of size $n_s$ instead of a large covariance matrix. However, the number of realizations required to provide reasonable solutions is very high and this results in higher storage and computational costs. In CEKF, the cross-covariance $\Skk H^T$ is propagated between iterations and as a result, both the storage and computational costs for CEKF scale as $\bigO(n_mn_s)$. The Fast Kalman Filter proposed in this paper also enjoys the same storage and computational costs as CEKF proposed in~\cite{li2014kalman}.

Figure~\ref{fig:errkalman} compares the error between the proposed fast Kalman filter (FKF) against the standard Kalman filter (KF) and the Ensemble Kalman Filter (EnKF) on a grid size $59\times 55$. It can be seen that the error of the FKF is comparable with the error of the KF. If all the eigenvalues were computed accurately this error would be close to machine precision. However, the randomized algorithm trades computational efficiency for accuracy (a further description is available in Appendix~\ref{sec:rand}). On the other hand, the Ensemble Kalman filter has low accuracy compared to the KF and the FKF even with $1000$ ensembles. This is because the matrix $\prior$ does not have rapidly decaying eigenvalues and therefore a larger ensembles are necessary to produce an accurate representation. The results are only compared on the coarsest grid because computing the ensembles can be expensive\footnote{Computing ensembles are implemented in this paper using dense Cholesky factorization which is prohibitively expensive on finer grids. However, a more computationally efficient approach would be to use the methods described in section~\ref{sec:sampling}.}. This is another shortcoming of the Ensemble based approach.

\begin{figure}[!ht]
\centering
\includegraphics[scale=0.5]{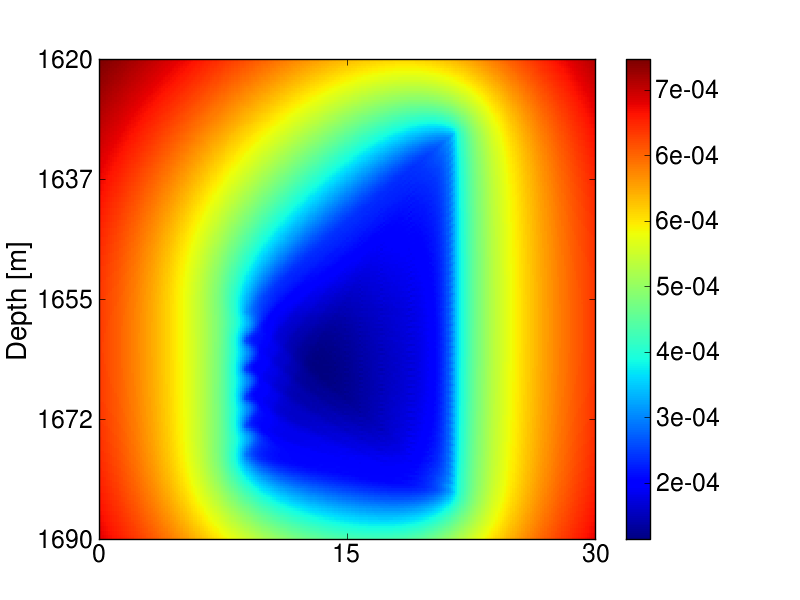}
\caption{Variance [TL$^{-1}$] of the computed solution at time $30$ hours after injection computed on the grid of size $234\times 217$.}
\label{fig:variance}
\end{figure}

\begin{figure}[!ht]
\centering
\includegraphics[scale=0.5]{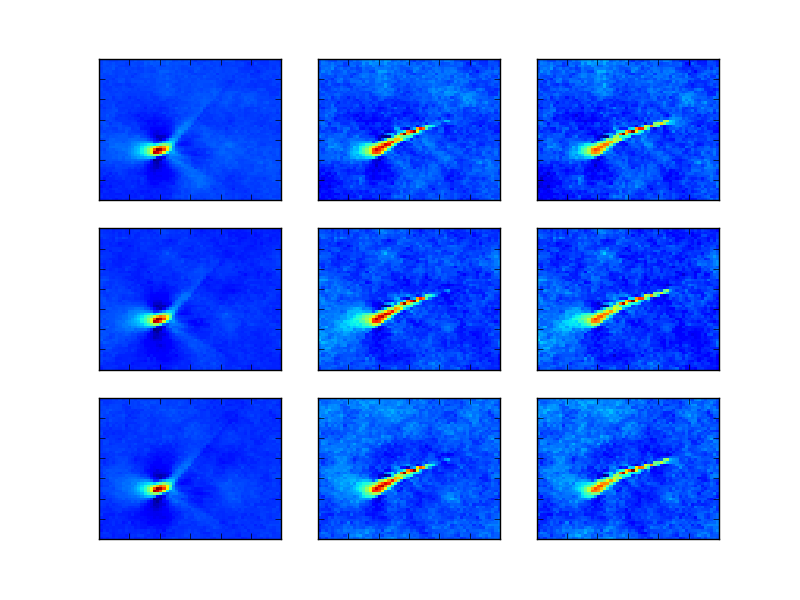}
\caption{Samples from the posterior distribution representing the seismic slowness [TL$^{-1}$]  of CO$_2$ plumes at $3$, $30$ and $60$ hours respectively computed on the grid of size $59\times 55$}
\label{fig:samples}
\end{figure}

We now address the issue of uncertainty quantification using the Kalman filter. In figure~\ref{fig:variance}, the estimated variance is plotted as a function of position, 30 hours after injection. The variance is computed as the diagonals of the covariance matrix $\Skk$. The variance is higher in regions that are not under the ray coverage. Another way to visualize the posterior covariance matrix $\Skk$ to understand the uncertainty and variability, is to visualize realizations from the distribution. In figure~\ref{fig:samples}, we plot samples from the distribution at $3,30$ and $60$ hours. The MVPs $\prior^{1/2}x$ and $\prior^{-1/2}x$ are computed by forming the Cholesky decomposition and as a results are limited to the coarsest grid of size $59\times 55$. A more scalable implementation can be obtained by using a contour integral based approach~\cite{hale2008computing}. 

\subsection{Extended Kalman Filter}\label{sec:extkalman}
In the previous subsection, we assumed that the measurement operator is time-invariant, i.e. $H_k=H$ for all time steps $k=1,\dots,N_t$. In terms of the application discussed, this is a consequence of the fact that the sources and receiver locations do not change with time. In order to demonstrate that our algorithm can be extended to the more general setting, i.e. $H_k$ changes with time, we consider a nonlinear transformation with repeated linearizations of the state, resulting in different measurement operators at each time step. We consider the following class of nonlinear transformations,  
\begin{equation}\label{eqn:nonlinear} \hat{s} = \alpha(s^{1/\alpha} - 1)  \qquad s = \left(\frac{\hat{s}+\alpha}{\alpha}\right)^{\alpha}\end{equation}
where $\alpha$ is a positive parameter that controls the degree of non-linearity of the transformation. The above transformation is known as Box-Cox transformation and is a useful mathematical tool for enforcing the non-negativity constraint. Another frequently used transformation is $\hat{s} = \log(s)$. In fact, in the limit $\alpha \rightarrow \infty$, we have that  $ \alpha(s^{1/\alpha} - 1) \rightarrow \log(s)$. The logarithm transformation encounters difficulties when the parameter $s$ is zero. As a result we prefer the power transformation with high values of $\alpha$ since $s$ is constrained to be non-negative and the mapping between $s$ and $\hat{s}$ is one-to-one. In the context of our application, since we are estimating the concentration of CO$_2$, this nonlinear transformation ensures that the reconstruction of the concentration remains positive leading to a more physically ``realistic'' situation.  

We assume that the CO$_2$ concentration $s_k$ and the seismic travel time measurements $y_k$ satisfy the following dynamical system 
\begin{align*}
s_{k+1} = & \quad F_ks_{k} + w_k  &\qquad w_k \sim {\cal N} (0,\prior)  \\ \nonumber
y_{k+1} = & \quad h(s_{k+1}) + v_k  &\qquad v_k \sim {\cal N} (0,\noise).  \\ \nonumber
\end{align*}
As before, we make the random-walk assumption i.e., $F_k = I$. Since the operator is nonlinear whenever $\alpha \neq 1$, the Kalman Filter approach cannot be used directly and we need to adopt an Extended Kalman Filter approach. Essentially, the nonlinear measurement operator is linearized about the current estimate using Taylor series expansion. We have 
\[ h(s) = h(s_k) + H_k(s-s_k) + \bigO(\normtwo{s-s_k}^2)\qquad H_k = \at{\frac{\partial h}{\partial s}}{s_k} = H\text{diag}\left(\frac{\partial s}{\partial \hat{s}}\right)\] 
As a result, we make the following modification in the update of the state vector $\skkpp$ and the posterior covariance matrix $\Skkpp$. The final algorithm is summarized in algorithm~\ref{alg:extkalman}.
  
\begin{align*}
\Skkpp \quad = & \quad \left[ \Skkp^{-1} + H_k^T\noise^{-1} H_k \right]^{-1} \\ 
\skkpp \quad = & \quad \skk+\Skkp H^T_k (H_k\Skk H_k+\noise)^{-1}\left(y_{k+1} - h(\skk)\right)  
\end{align*}

\begin{algorithm}
\begin{algorithmic}[1]
\REQUIRE Measurement operators $H_k$ and measurements $y_k$ for $k=1,\dots,N_t$, system noise covariance $\prior$ and measurement noise $\noise$ 
\FOR {$k=1,\dots,N_t$}
\STATE $\alpha_{k+1}=\alpha_k+1$ and $\bar{D}_k = \alpha_{k+1}^{-1}( \alpha_{k+1}I - D_k)^{-1}D_k$
\STATE Compute the Jacobian $H_k = \at{\frac{\partial h}{\partial s}}{s_k} = H\text{diag}\left(\at{\frac{\partial s}{\partial \hat{s}}}{\skk}\right)$
\STATE Compute the generalized eigendecomposition \[ H_k^T\noise^{-1}H_k = \prior^{-1}U_k \Lambda_k U_k^T\prior^{-1}\qquad \text{with}\qquad  U_k^T\prior^{-1}U_k = I\]
\STATE Compute $F_k\define \Skk H_k^T  = \alpha_{k+1}\prior H_k^T - W_kD_k(W_k^TH_k)$ 
\STATE Compute $\skkp =  \skk + F_k(H_kF_k+\noise)^{-1}\left(y_k-h(\skk)\right)$. \\
\COMMENT {//If necessary, do additional linearization steps.}
\STATE $[W_{k+1},\hat{D}] = \text{AddLowRank}\left(W_k,\bar{D}_k, U_k,\Lambda_k,\prior^{-1}\right)$
\STATE $D_{k+1} \define \alpha_{k+1} (I + \alpha_{k+1}^{-1}\hat{D}^{-1})^{-1} $
\ENDFOR
\end{algorithmic}
\caption{Fast Extended Kalman Filter for random walk forecast model}
\label{alg:extkalman}
\end{algorithm}

\begin{figure}[!ht]
\centering
\includegraphics[scale=0.29]{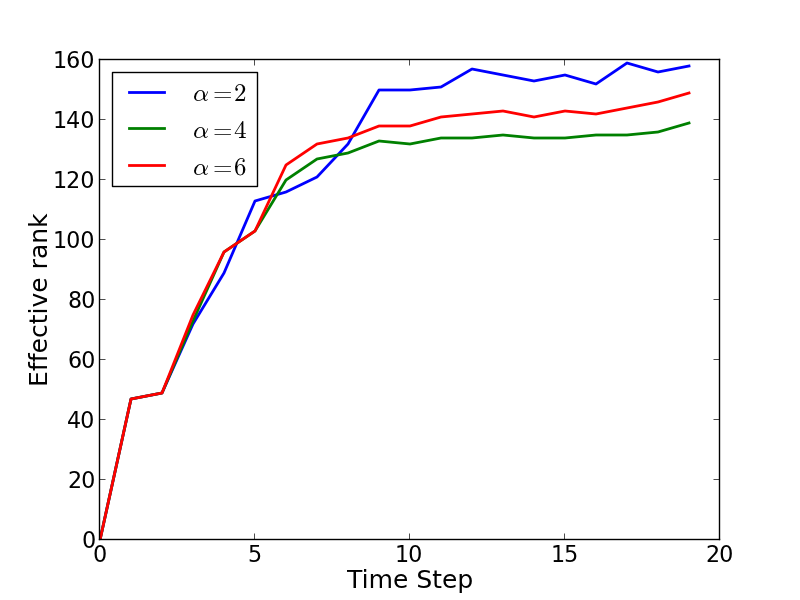}
\includegraphics[scale=0.29]{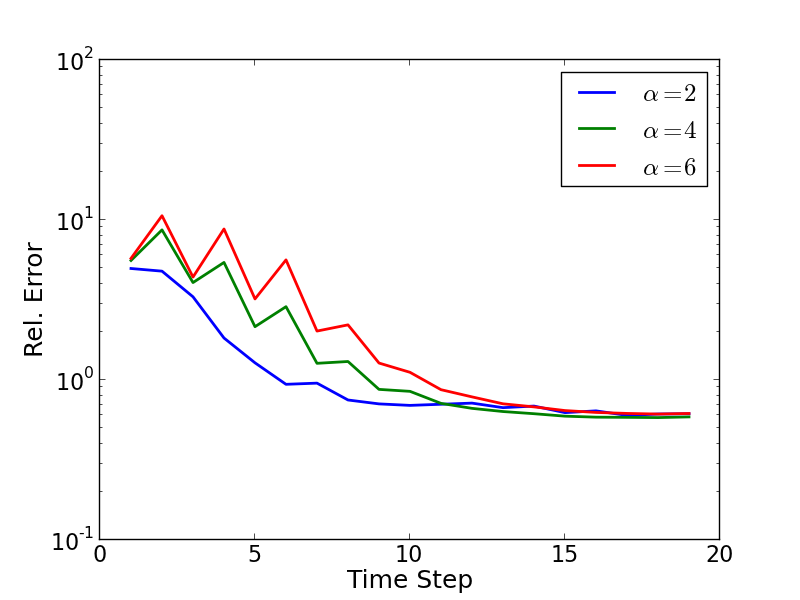}
\caption{(left) Effective rank of the matrix $W_k$ as a function of the time step $k$ for different values of $\alpha$ (the parameter of the nonlinearity transformation). A relative tolerance of $10^{-5}$ was used to truncate the effective rank at each time step. (right) the relative $L^2$ error in the extended Kalman filter compared to the ``true'' solution for multiple values of $\alpha$.}
\label{fig:extkalman}
\end{figure}

The results of the algorithm~\ref{alg:extkalman} on the same application described in section~\ref{sec:setup} is presented in figure~\ref{fig:extkalman}. The nonlinear transformation described in equation~\eqref{eqn:nonlinear} was used with three different values of $\alpha =2,4,6$ where increasing values of $\alpha$ correspond to higher non-linearity. A covariance kernel was used  to be $\kappa(r) = \theta\exp\left(-\frac{r^p}{l^p} \right)$, with $p = 1$ and $\theta=10^{-5}$. We assume that $\noise = \sigma^2I$ with $\sigma^2 = 2\times 10^{-4}$. The kernel corresponding to $p=1$ is smoother than $p=1/2$ (used in the previous section) and we expect the eigenvalues of the generalized eigenvalue problem~\eqref{eqn:ghep} to decay more rapidly, resulting in a more efficient representation of the state covariance matrix $\Skk$. A relative tolerance of $10^{-5}$ was used to truncate the effective rank of $W_k$. The full rank of the perturbative term grows linearly as $n_m\times t$, where $t$ is the number of time-steps elapsed and $n_m=288$ is the number of measurements per time-step. As can be seen the effective rank of the low-rank perturbation increases steadily and then reaches a plateau because no additional information in terms of eigenmodes comes through the combination of the measurement operator $H_k$ and the system noise covariance $\prior$. Further, the error in the first few linearization steps is high because the linearization about the  initial field (assumed to be zero everywhere) is erroneous and the error decrease when additional information enters the assimilation algorithm. Furthermore, the error is comparable to the Kalman Filter which has already been reported in figure~\ref{fig:errkalman}.

We conclude this section with a discussion on the computational cost of the proposed fast algorithm for the Extended Kalman filter. Several components are necessary for efficient application of algorithm~\ref{alg:extkalman} 1) the computational cost of the generalized eigenvalue problem at every iteration in equation~\eqref{eqn:ghep} needs to be scalable, 2) the rank of the low-rank perturbative matrix needs to be bounded and 3) the efficient truncation of the low-rank addition described in~\ref{alg:lowrankaddition}. We have already demonstrated the scalability of the eigenvalue calculations and the bounded effective rank of the low-rank perturbation. In summary, the computational cost of updating a step is $\bigO(n_s\log n_s + r^2n_s)$, where $r$ is the rank of the perturbation. However, the algorithm~\ref{alg:lowrankaddition} requires repeated application involving the inverse of $\prior$. This computation is also scalable as $\bigO(n_s\log n_s)$ because of the arguments made in section~\ref{sec:cov}; however, the pre-factor in front of the computation is extremely high. Further research is necessary to develop an implementation with a lower pre-factor making it competitive. This issue will be further explored in a forthcoming paper~\cite{saibaba2014fast} and will not be discussed further here.  

\section{Conclusions}
We have presented a fast algorithm for updating the estimate of system state and the associated uncertainty for an application arising from monitoring CO$_2$ plumes in the subsurface using time-lapse seismic signals. The key step was an efficient representation of the state covariance matrix as a low-rank perturbation of the system noise covariance matrix which is appropriately weighted. The low-rank perturbation combines information about the measurement operator and the system noise covariance matrix by solving a generalized eigenvalue problem that takes advantage of the eigenvalue decay in one or both matrices. When the measurement operator is time-invariant, the weights associated with the representation can be efficiently updated at negligible cost since additional information only propagates through the dominant eigenmodes of the generalized eigenvalue problem, which have been precomputed. When the measurement operator changes in time, the additional information can be incorporated by computing the generalized eigendecomposition at each time step. Consequently, the rank of the low-rank perturbation grows linearly. In order to alleviate the  computational burden we have proposed a method to truncate the rank of the perturbation. The resulting algorithm (although expensive) still scales almost linearly with the number of state variables. The additional advantage of this efficient representation of the state covariance matrix is the ability to efficiently compute various uncertainty measures which are scalar functions of the state covariance matrix. These have been discussed in section~\ref{sec:uncert}. We have demonstrated the scalability and accuracy of our algorithm by means of numerical examples.

In future work, we would like to consider an extension of our ideas to dynamical systems for which the state transition matrix $F_k$ is not equal to identity, for time steps $k=1,\dots,N_t$. For example, our approach can be extended to systems for which the state transition matrices $F_k$ are normal and commute with $\prior$ as would be the case where these matrices are (block) circulant/Toeplitz. An example for which circulant matrices are relevant, are a constant coefficient PDE, for e.g., the advection-diffusion equation with periodic boundary conditions. Another possible avenue is the use of Hierarchical matrix arithmetic to approximate expensive matrix-matrix products (such as those that arise in the prediction step $\Sigma_{k+1|k} = F_k\Sigma_{k|k} F_k^T + \prior$) which can be done in almost linear time, when all the relevant matrices are in the ${\mathcal{H}}$-matrix of ${\mathcal{H}}^2$-matrix format. Examples include dynamical systems for which the governing equations can be described using parabolic PDEs. In this example, the matrices that arise from the discretized operators (which are elliptic) have been shown to be efficiently approximated in the $\mathcal{H}$-matrix format~\cite{Grasedyck03constructionand}. Other dynamical systems for which our analysis is relevant has been summarized in~\cite{pnevmatikakis2013fast}. We note that the Ensemble Kalman filter may also work well for these problems and future work could compare the strengths and weakness of both approaches. The extension of these ideas to other data assimilation methods such as 4D-Var data assimilation would also be interesting to consider.

\section{Acknowledgments}
The authors would also like to thank Dr. Jonathan B. Ajo-Franklin, Thomas M. Daley, and Christine Doughty from the Lawrence Berkeley Lab for sharing TOUGH2 and rock physics simulation data. The first author would also like to thank Judith Y. Li for her help with the data set.  We would also like to thank the anonymous reviewers for their careful reading of the manuscript and their comments that helped improve the presentation of the paper. Subsection~\ref{sec:parallel} has been added as per their suggestion.

\appendix
\section{Adding low-rank matrices}\label{sec:addlowrank}

In this appendix we derive an efficient algorithm for adding low rank matrices such that the resulting low rank representation is a generalized eigendecomposition. This will be useful in section~\ref{sec:nonconstantH}. Consider two low rank matrices $UD_UU^T$ and $VD_VV^T$ where $U$ and $V$ satisfy $U^TBU=I$ and $V^TBV=I$ and $D_U$ and $D_V$ are diagonal matrices. We want to represent the result of adding two low rank matrices $A = UD_UU^T + VD_VV^T $ as $A =  WD_WW^T$, where $W^TBW = I$. In general, adding two low rank matrices will produce a low rank matrix whose rank is the less than equal to the sum of the ranks of low rank matrices. In order to make the representation more efficient, we consider truncating the singular values that are below a certain threshold. The algorithm is summarized in~\ref{alg:lowrankaddition}.

\begin{algorithm}[!ht]
\begin{algorithmic}[1]
\REQUIRE Low rank matrices $UD_UU^T$ and $VD_VV^T$, an symmetric positive definite matrix $B$ and a tolerance $\text{tol}$ \\ 
\COMMENT {// Assume that $U^TBU=I$ and $V^TBV=I$.}
\STATE Compute $ V_t \leftarrow V - UU^TBV $ \COMMENT{//Block Gram-Schmidt}
\STATE Compute the QR factorization $\hat{V}R = V_t$ such that $\hat{V}^TB\hat{V}=I$ \\
\COMMENT {// $W = [U,\hat{V}]$ forms a B-orthonormal basis for $\Span{[U,V]}$.}
\STATE  Form the matrix \[  M =  \begin{bmatrix} I \\ 0 \end{bmatrix} D_U \begin{bmatrix} I & 0 \end{bmatrix} + \begin{bmatrix} U^TBV \\ \hat{V}^TBV\end{bmatrix} D_V \begin{bmatrix} V^TBU & {V}^TB\hat{V}\end{bmatrix}  \]
\STATE Compute eigendecomposition $M = S\Lambda S^T$. Truncate eigenvalues below the tolerance.
\STATE $W\leftarrow WS$ and $D_W =\Lambda$
\end{algorithmic}
\caption{Adding low rank matrices, $[W,D_W] = \text{AddLowRank}(U,D_U,V,D_V,B,\text{tol})$}
\label{alg:lowrankaddition}
\end{algorithm}

\section{Appendix: Solving the GHEP}\label{sec:rand}
In section~\ref{sec:fastkalman} we needed to repeatedly solve the generalized Hermitian eigenvalue problem (GHEP) $H_k\noise^{-1}H_kx = \lambda \prior^{-1}x$ in order to compute the dominant eigenmodes. This problem has been previously developed in the following references~\cite{saibaba2013randomized,saibaba2013uncertainty,saibaba2013fast}. For the sake of completion, we reproduce the section on solving the GHEP from~\cite[chapter 6, section 6.5]{saibaba2013fast}.

We briefly review the randomized algorithm described in~\cite{saibaba2013randomized} for computing dominant eigenmodes of the GHEP $Ax = \lambda Bx$. In the context of solving the problem~\eqref{eqn:ghep}, we have $A\define H^T\noise^{-1}H$ and $B \define \prior^{-1}$. The key observation is that the matrix $C \define B^{-1}A$ is symmetric with respect to the $B$-inner product $\langle x,y\rangle_B = y^TBx$. Suppose we wanted to compute the $k$ largest generalized eigenpairs of $Ax=\lambda B x$. We assume that $B^{-1}x$ is easier to compute than $Bx$. This is certainly the case since $B=\prior^{-1}$. The randomized algorithm~\ref{alg:randomizedghepuncert} calculates a matrix $Q$, which is $B$-orthonormal and approximately spans the column space of $C$, i.e. satisfies the following error bound $\normB{(I-QQ^*B)C} \leq \varepsilon $. Given such a matrix $Q$, it can be shown that $\normB{ A\approx (BQ) (Q^*AQ) (BQ)^*} \leq 2\varepsilon$, i.e. $A\approx (BQ) (Q^*AQ) (BQ)^*$. As a result, a symmetric rank-$k$ approximation can be computed, from which the approximate eigendecomposition can be computed.

To produce a symmetric rank-$k$ approximation, the algorithm proceeds as follows: first, we sample a matrix with entries randomly chosen from $\normal (0,1)$, $\Omega \in \mathbb{R}^{n\times r}$. We choose $r = k +p$, where $p$ is an oversampling factor, which we choose to be $20$. Form $\bar{Y} = A\Omega$. Then, we compute the QR factorization of  $\bar{Y} =\bar{Q}R$ such that $\bar{Q}^*B^{-1}\bar{Q} = I$. This can be accomplished by modified Gram-Schmidt algorithm with $\langle\cdot,\cdot\rangle_{B^{-1}}$ inner products. Then, compute $Q\define B^{-1}\bar{Q}$ which is now $B$-orthonormal. Then, we form $T\define Q^*AQ$ and compute its eigenvalue decomposition $T = S\Lambda S^*$. We then have the approximate generalized eigendecomposition 
\[ A \approx U\Lambda U^* \qquad U = \bar{Q}S\]
Here, $U$ is also $B$-orthonormal. The cost for computing the $k$-largest modes just involves $2r$ MVPs with $A$ and $2r$ MVPs with $B^{-1}$ and an additional cost that is $\bigO(r^2n)$. The cost of a second round of MVPs with $A$ while computing $T$ can be avoided using the following observation:
\[ \Omega^*\bar{Y} = \Omega^*A\Omega \approx (\Omega^* BQ) T (BQ^*\Omega)\]
Therefore, $T$ can be computed as $T \approx (\Omega^*\bar{Q})^{-1}(\Omega^*\bar{Y}) (\bar{Q}^*\Omega)^{-1}$. This is summarized in algorithm~\ref{alg:randomizedghepuncert}.

\begin{algorithm}[!ht]
\begin{algorithmic}[1]
 \REQUIRE  matrices $A$, $B$ and $\Omega \in \mathbb{R}^{n\times (k+p)}$ is a Gaussian random matrix. Here $A,B \in \mathbb{C}^{n\times n}$, k is the desired rank, $p\sim 20$ is an oversampling factor. 
\STATE  Compute $\bar{Y} = A\Omega$
\STATE  Form QR factorization $\bar{Y}=\hat{Q}R$ such that $\hat{Q}^*B^{-1}\hat{Q} = I$ using algorithm~\cite[algorithm 2]{saibaba2013randomized} (with $W=B^{-1}$).
\STATE Compute $Q = B^{-1}\hat{Q}$ so that $Q^*B Q = I$.
\STATE  Form $ T = Q^*A Q$ or $\approx (\Omega^*\bar{Q})^{-1}(\Omega^*\bar{Y}) (\bar{Q}^*\Omega)^{-1}$,  and
\STATE  Compute the eigenvalue decomposition $T = S\Lambda S^*$. Keep the $k$ largest eigenmodes.
\STATE  \textbf{Return:} Matrices $U \in \mathbb{R}^{n\times k}$ and $\Lambda\in \mathbb{R}^{k\times k}$ 
 that satisfy
\[ A  \approx  (BU) \Lambda (BU)^*\qquad\text{with}\qquad U = QS \]
\end{algorithmic}
\caption{Randomized algorithm for GHEP }
\label{alg:randomizedghepuncert}
\end{algorithm}

The efficiency and accuracy of this algorithm has been studied for several test problems and the reader is referred to~\cite{saibaba2013randomized}. Here we summarize the main conclusions. The error in the low-rank approximation is $$\normB{(I-QQ^*B)B^{-1}A} \leq c\normtwo{B^{-1}}\sigma_{B,k+1}(B^{-1}A)$$
where $c$ is a constant that depends on $n,k$ and $p$ and is independent of the spectrum of the matrices. $\sigma_{B,k+1}$ is the $(k+1)$-th generalized singular value of the matrix $B^{-1}A$. Since it is hard to compute the generalized singular values, an randomized estimator for the error in the low-rank representation is also proposed and analyzed in~\cite{saibaba2013randomized}. Given error in the low-rank representation $\normB{(I-QQ^*B)B^{-1}A}$, it can be shown that the error in approximating the true eigenvalue and eigenvector satisfies the following error bounds  
\[|\lambda-\tilde{\lambda}| \leq \min\{2\varepsilon,\frac{4\varepsilon^2}{\delta}\}\qquad \sin\angle_B (u,\tilde{u}) \leq \frac{2\varepsilon}{\delta}\]

where $\delta = \min_{\lambda_i\neq\lambda} |\tilde{\lambda}-\lambda_i|$ is the gap between the approximate eigenvalue $\tilde{\lambda}$ and any other eigenvalue and $\angle_B(x,y) = \arccos \frac{|\langle x,y\rangle_B|}{\normB{x}\normB{y}}$. This result states that the accuracy in the eigenvalue/eigenvector calculations depends not only on the accuracy of the low-rank representations but also on the spectral gap $\delta$. When the eigenvalues are clustered, the spectral gap is small and the eigenvalue calculations are accurate as long as the error in the low-rank representation is small. However, in this case the resulting eigenvector calculations maybe inaccurate because the parameter $\delta$ appears in the denominator for the approximation of the angle between the true and approximate eigenvector.

\bibliographystyle{plain}
\bibliography{kalman}

\end{document}